\documentclass{amsart}
\usepackage{amssymb}
\usepackage{amsmath}
\usepackage{amsfonts}

\newtheorem{theorem}{Theorem}
\theoremstyle{plain}
\newtheorem{acknowledgement}{Acknowledgement}

\newtheorem{corollary}{Corollary}

\newtheorem{definition}{Definition}
\newtheorem{example}{Example}

\newtheorem{lemma}{Lemma}

\newtheorem{remark}{Remark}

\numberwithin{equation}{section}
\input{tcilatex}

\begin{document}
\title[LDP for Enhanced Gaussian Processes]{Large Deviation Principle for
Enhanced Gaussian Processes}
\author{Peter Friz, Nicolas Victoir}
\email{P.K.Friz@statslab.cam.ac.uk}

\begin{abstract}
We study large deviation principles for Gaussian processes lifted to the
free nilpotent group of step N. We apply this to a large class of Gaussian
processes lifted to geometric rough paths. A large deviation principle for
enhanced (fractional) Brownian motion, in H\"{o}lder- or modulus topology,
appears as special case.
\end{abstract}

\maketitle

\textrm{R\'{e}sum\'{e}: }\textit{Nous etudions les principes de grandes
deviations pour les "rough paths Gaussiens", pour une large classe de
processus Gaussiens. Cette classe contient le mouvement brownien
fractionnaire.}

\section{Introduction}

We prove large deviation principles for $d$-dimensional Gaussian processes
lifted to the free nilpotent group of step $N$. The example we have in mind
is a class of zero mean Gaussian processes subject to certain conditions on
the covariance, \cite{CQ, LQ}.

After recalls on nilpotent groups and Wiener It\^{o} chaos, we give simple
conditions under which a large deviation principle holds. To illustrate the
method, we quickly check these conditions for enhanced Brownian motion (see 
\cite{LQZ, FV03} for earlier approaches). We then apply our methodology to
the class of enhanced Gaussian processes considered in \cite{CQ, LQ}.
Enhanced fractional Brownian motion appears as a special case and we
strengthen the results in \cite{MS, CFV}.

\subsection{The Free Nilpotent Group and Rough Paths}

Let $N\geq 0$. The truncated tensor algebra of degree $N$ is given by the
direct sum 
\begin{equation*}
T^{N}\left( \mathbb{R}^{d}\right) =\mathbb{R}\oplus \mathbb{R}^{d}\oplus
...\oplus \left( \mathbb{R}^{d}\right) ^{\otimes N}.
\end{equation*}%
With the usual scalar product, vector addition, and tensor product $\otimes $
the space $T^{N}\left( \mathbb{R}^{d}\right) $ is an algebra. Let $\pi _{i}$
denote the canonical projection from $T^{N}\left( \mathbb{R}^{d}\right) $
onto $\left( \mathbb{R}^{d}\right) ^{\otimes i}$.

Let $q\in \lbrack 1,2)$ and $x\in $ $C_{0}^{q\text{-var}}\left( \left[ 0,1%
\right] ,\mathbb{R}^{d}\right) $, the space of continuous $\mathbb{R}^{d}$%
-valued paths of bounded $q$-variation started at the origin. The lift $%
S_{N}(x):[0,1]\rightarrow T^{N}\left( \mathbb{R}^{d}\right) $ is defined as 
\begin{equation*}
S_{N}(x)_{t}=1+\sum_{i=1}^{N}\int_{0<s_{1}<...<s_{i}<t}dx_{s_{1}}\otimes
...\otimes dx_{s_{i}},
\end{equation*}%
the iterated integrals are Young integrals. Then 
\begin{equation*}
G^{N}\left( \mathbb{R}^{d}\right) =\left\{ g\in T^{N}\left( \mathbb{R}%
^{d}\right) :\exists x\in C_{0}^{1\text{-var}}\left( \left[ 0,1\right] ,%
\mathbb{R}^{d}\right) :\text{ }g=S_{N}(x)_{1}\text{ }\right\}
\end{equation*}%
is a submanifold of $T^{N}\left( \mathbb{R}^{d}\right) $ and, in fact, a Lie
group with product $\otimes $, called the free nilpotent group of step $N$.
The dilation operator $\delta :\mathbb{R}\times G^{N}\left( \mathbb{R}%
^{d}\right) \rightarrow G^{N}\left( \mathbb{R}^{d}\right) $ is defined by 
\begin{equation*}
\pi _{i}\left( \delta _{\lambda }(g)\right) =\lambda ^{i}\pi
_{i}(g),\,\,\,i=0,...,N.
\end{equation*}%
and a continuous norm on $G^{N}\left( \mathbb{R}^{d}\right) $, homogenous
with respect to $\delta $, is given 
\begin{equation*}
\left\Vert g\right\Vert =\inf \left\{ \text{length}(x):x\in C_{0}^{1\text{%
-var}}\left( \left[ 0,1\right] ,\mathbb{R}^{d}\right)
,S_{N}(x)_{1}=g\right\} .
\end{equation*}%
By equivalence of such norms there exists a constant $K_{N}$ such that%
\begin{equation}
\frac{1}{K_{N}}\max_{i=1,...,N}|\pi _{i}(g)|^{1/i}\leq \left\Vert
g\right\Vert \leq K_{N}\max_{i=1,...,N}|\pi _{i}(g)|^{1/i}.  \label{eqNorms}
\end{equation}

The norm $\left\Vert \cdot \right\Vert $ induces a metric on $G^{N}\left( 
\mathbb{R}^{d}\right) $ known as Carnot-Caratheodory metric,%
\begin{equation*}
d(g,h)=\left\Vert g^{-1}\otimes h\right\Vert .
\end{equation*}

Let $x,y\in $ $C_{0}\left( [0,1],G^{N}\left( \mathbb{R}^{d}\right) \right) $%
, the space of continuous $G^{N}\left( \mathbb{R}^{d}\right) $-valued paths
started at the neutral element of $\left( G^{N}\left( \mathbb{R}^{d}\right)
,\otimes \right) $. We define a supremum distance, 
\begin{equation*}
d_{\infty }(x,y)=\sup_{t\in \lbrack 0,1]}d(x_{t},y_{t}),
\end{equation*}%
a H\"{o}lder distance, $\alpha \in \left[ 0,1\right] ,$ 
\begin{equation*}
d_{\alpha -\text{H\"{o}lder}}(x,y)=\sup_{0\leq s<t\leq 1}\frac{%
d(x_{s}^{-1}\otimes x_{t},y_{s}^{-1}\otimes y_{t})}{|t-s|^{\alpha }},
\end{equation*}%
write 
\begin{equation*}
\left\Vert x\right\Vert _{\alpha -\text{H\"{o}lder}}=\sup_{0\leq s<t\leq 1}%
\frac{\left\Vert x_{s}^{-1}\otimes x_{t}\right\Vert }{|t-s|^{\alpha }}%
=\sup_{0\leq s<t\leq 1}\frac{d\left( x_{s},x_{t}\right) }{|t-s|^{\alpha }},
\end{equation*}

and $Lip$ instead of $1$-H\"{o}lder. Similarly, one defines $p$-variation
regularity and distance for $G^{N}\left( \mathbb{R}^{d}\right) $-valued
paths. We refer to \cite{FV04} for a more detailed discussion of these
topics.

The interest for all this comes from T. Lyons' rough path theory, \cite%
{Ly,LQ}, a deterministic theory of control differential equations ("Rough
Differential Equation"), driven by a continuous $G^{[p]}\left( \mathbb{R}%
^{d}\right) $-valued path of finite $p$-variation. Brownian motion and L\'{e}%
vy's area can be viewed as a $G^{[p]}\left( \mathbb{R}^{d}\right) $-valued
path of finite $p$-variation with $p>2$ and the corresponding RDE solution
turns out to be a classical Stratonovich SDE solution. This gives a general
recipe how to make sense of differential equations driven by an \textit{%
arbitrary} stochastic process $X$ with continuous sample path of finite $p$%
-variation: find an appropriate lift of $X$ to a $G^{[p]}\left( \mathbb{R}%
^{d}\right) $-valued path $\mathbf{X}$ with finite $p$-variation use Lyons'
machinery.

\subsection{Gaussian Process}

The Wiener space $C_{0}\left( \left[ 0,1\right] ,\mathbb{R}^{d}\right) $ is
the space of continuous function started at $0$; \ the sup norm induces a
topology and hence a Borel $\sigma $-algebra. We assume there is a
probability measure $\mathbb{P}$ such that the coordinate process $X$,
defined by $X\left( \omega \right) _{t}=\omega \left( t\right) $, is a
centered Gaussian process with covariance $c\left( s,t\right) .$ We also
define $\mathcal{H}$ to be the reproducing kernel Hilbert space associated
to $\mathbb{P}$, i.e. the closure of the set of the linear span of the
functions $s\rightarrow c\left( t_{i},s\right) ,$ under the Hilbert product 
\begin{equation*}
\left\langle c\left( t_{i},.\right) ,c\left( t_{j},.\right) \right\rangle
=c\left( t_{i},t_{j}\right) .
\end{equation*}%
Let $\emph{i}$ denote the canonical injection from the reproducing kernel
Hilbert space to the Wiener space. We define $\mathcal{C}_{n}(B)$ the $B$%
-valued homogeneous Wiener chaos of degree $n$, where $(B,|.|_{B})$ is a
given Banach space. We also define $C_{n}(B)$ the sum of the first $n$
homogeneous Wiener chaos. We refer to \cite{Ledoux, B1, B2, B3} for more
details on Gaussian spaces and Wiener chaos.

\subsection{Enhanced Gaussian Process}

We make the assumption that there exists a lift of the Gaussian process $X$
to a process $\mathbf{X}$ with values in $G^{N}\left( \mathbb{R}^{d}\right) $
for some fixed $N\geq 1$. Such lifts have been constructed first by Coutin
and Qian \cite{CQ} as a.s. limits in $p$-variation of (canonically lifted)
dyadic piecewise linear approximations of $X$. We shall be less specific
here and extract those properites of the approximations that we need in the
sequel. We make the following

\textbf{Assumption:} There exists a sequence of continuous linear maps $%
\left\{ \Phi _{m}\right\} $ from $C\left( \left[ 0,1\right] ,\mathbb{R}%
^{d}\right) $ (with uniform topology) onto the space of bounded variation
path from $\left[ 0,1\right] $ into $\mathbb{R}^{d}$ (with $1$-variation
topology) such that

\begin{enumerate}
\item $\lim_{m\rightarrow \infty }\Phi _{m}\left( x\right) =x$ for all
continuous paths $x$ with respect to uniform topology;

\item the uniform distance between $S_{N}\circ \Phi _{m}\left( X\right) $
and $\mathbf{X}$ converges to $0$ in probability.
\end{enumerate}

In fact, condition 2 is equivalent to the seemingly stronger condition of $%
L^{p}$-convergence.

\begin{lemma}
\label{proba_>Lq}Convergence of $d_{\infty }\left( S_{N}\circ \Phi
_{m}\left( X\right) ,\mathbf{X}\right) $ to zero in probability is
equivalent to convergence in in $L^{p}$ for all $p\in \lbrack 1,\infty )$.
\end{lemma}

\begin{proof}
Theorem III.2 in \cite{S} and equation (\ref{eqNorms}).
\end{proof}

\begin{remark}
Assumption 2 is always satisfied when $N=1$. In the applications discussed
in later sections, $N$ will be related to the regularity of $X$.
\end{remark}

\begin{example}
Let $\left\{ \Phi _{m}\right\} $ be the piecewise linear approximations
based on a sequence of dissections $\left\{ D_{m}\right\} $ with mesh $%
\left\vert D_{m}\right\vert \rightarrow 0.$ If $X$ denotes fractional
Brownian motion of parameter $H>\frac{1}{4},$ then the lift $\mathbf{X}$ for 
$N=[1/H]$ and was first constructed in \cite{CQ}, see also \cite{CFV}.
\end{example}

\begin{remark}
As is well known (e.g. Theorem 1 in \cite{JaKa}) any continuous Gaussian
process $X$ on $\left[ 0,1\right] $ is the uniform limit of%
\begin{equation*}
X_{t}^{\left( m\right) }\equiv \sum_{j=1}^{m}\xi _{j}\left( \omega \right)
\psi _{j}\left( t\right)
\end{equation*}%
where $\left\{ \psi _{j}\right\} $ is an orthonormal basis for the
reproducing Kernel Hilbert space $\mathcal{H}$ and the $\left\{ \xi
_{j}\right\} $ are i.i.d. standard Gaussian given by the image of $\psi _{j}$
under the isometric isomorphsim $\theta :\mathcal{H}\rightarrow \mathcal{L}%
_{2}\left( X_{t}:t\in \left[ 0,1\right] \right) $, the closure of $\left\{
X_{t}:\left[ 0,1\right] \right\} $ in $L^{2}\left( \mathbb{P}\right) $. The
Hilbert structure of $\mathcal{H}$ implies that for all $h\in \mathcal{H}$, 
\begin{equation*}
h^{\left( m\right) }\equiv \sum_{j=1}^{m}\left\langle \psi
_{j},h\right\rangle \psi _{j}\rightarrow h\text{ in }\mathcal{H}
\end{equation*}%
and hence uniformly. The construction of $\mathbf{X}$ can be based on such
(or similar) approximations\footnote{%
At least for, $N\leq 2$ it is shown in \cite{CV} that reproducing kernel -
and piecewise linear approximations yield the same lifted process.}, see 
\cite{CV,FV04, MS2, FD}, but they do not satisfy our assumption. It may be
possible to adapt the subsequent proofs, aimed to establish a large
deviation principle for $\mathbf{X}$, to such approximations. As this
requires further work without improving the results we shall not pursue this
further here.
\end{remark}

\section{Wiener Chaos and $G^{N}\left( \mathbb{R}^{d}\right) $-valued Paths 
\label{chaossection}}

The hypercontractivity of the Ornstein-Uhlenbeck semigroup leads to a useful

\begin{lemma}
\label{Hyper}Let $Z_{n}$ be a random variable in $\mathcal{C}_{n}(B)$.
Assume\ \ $1<p<q<\infty $. Then 
\begin{equation*}
\left\Vert Z_{n}\right\Vert _{L^{q}\left( \mathbb{P};B\right) }\leq \left( 
\frac{q-1}{p-1}\right) ^{\frac{n}{2}}\left\Vert Z_{n}\right\Vert
_{L^{p}\left( \mathbb{P};B\right) }.
\end{equation*}
\end{lemma}

\begin{proof}
Let $P_{t}$ denote the Ornstein-Uhlenbeck semigroup. Whenever $1<p<q<\infty $
and $t>0$ satisfies%
\begin{equation}
e^{t}\geq \left( \frac{q-1}{p-1}\right) ^{1/2},  \label{ChoiceOft}
\end{equation}%
then, for all $f\in L^{p}\left( \mathbb{P};B\right) $,%
\begin{equation*}
\left\Vert P_{t}f\right\Vert _{L^{q}}\leq \left\Vert f\right\Vert _{L^{p}}.
\end{equation*}%
We also recall that%
\begin{equation*}
P_{t}Z_{n}=e^{-nt}Z_{n}\text{.}
\end{equation*}%
See \cite{Ledoux} for this two results. We now choose $t$ such that equality
holds in (\ref{ChoiceOft}) and find%
\begin{equation*}
\left( \frac{q-1}{p-1}\right) ^{-n/2}\left\Vert Z_{n}\right\Vert
_{L^{q}}=e^{-nt}\left\Vert Z_{n}\right\Vert _{L^{q}}=\left\Vert
P_{t}Z_{n}\right\Vert _{L^{q}}\leq \left\Vert Z_{n}\right\Vert _{L^{p}}.
\end{equation*}
\end{proof}

It is well-known that $L^{p}\left( \mathbb{P};B\right) $- and $L^{q}\left( 
\mathbb{P};B\right) $-norms are equivalent on%
\begin{equation*}
C_{n}(B)=\mathcal{C}_{1}(B)+...+\mathcal{C}_{n}(B).
\end{equation*}%
The next lemma quantifies this equivalence.

\begin{lemma}
\label{LpLqEquivalence}Let $n\in \mathbb{N}$ and $Z$ be a random variable in 
$C_{n}(B)$. Assume $2\leq p\leq q<\infty $. Then there exists\ a constant $%
M_{n}$ such that%
\begin{equation*}
\left\Vert Z\right\Vert _{L^{p}\left( \mathbb{P};B\right) }\leq \left\Vert
Z\right\Vert _{L^{q}\left( \mathbb{P};B\right) }\leq M_{n}\left( q-1\right)
^{\frac{n}{2}}\left\Vert Z\right\Vert _{L^{p}\left( \mathbb{P};B\right) }
\end{equation*}
\end{lemma}

\begin{proof}
Only the second inequality requires a proof. Write $Z=\sum_{i=0}^{n}Z_{i}$
with $Z_{i}\in \mathcal{C}_{i}$. From Lemma \ref{Hyper}, for all $i\leq n$, 
\begin{equation}
\left\Vert Z_{i}\right\Vert _{L^{q}\left( \mathbb{P};B\right) }\leq \left( 
\frac{q-1}{p-1}\right) ^{i/2}\left\Vert Z_{i}\right\Vert _{L^{p}\left( 
\mathbb{P};B\right) }.  \label{hypercontractivity}
\end{equation}%
The sum $C_{n}(B)=\mathcal{C}_{1}(B)+...+\mathcal{C}_{n}(B)$ is topological
direct in $L^{0}\left( \mathbb{P};B\right) $, see \cite[page 6]{B3} for
instance, which implies that the projection $Z\mapsto $ $Z_{i}$ is
continuous in $L^{2}\left( \mathbb{P};B\right) $. It follows that $%
\left\Vert Z_{i}\right\Vert _{L^{2}\left( \mathbb{P};B\right) }\leq
c\left\Vert Z\right\Vert _{L^{2}\left( \mathbb{P};B\right) }$ for some
constant $c=c\left( n\right) $. Then, for $p\geq 2$,

\begin{equation}
\left\Vert Z_{i}\right\Vert _{L^{p}\left( \mathbb{P};B\right) }\leq \left(
p-1\right) ^{i/2}\left\Vert Z_{i}\right\Vert _{L^{2}\left( \mathbb{P}%
;B\right) }\leq c\left( p-1\right) ^{i/2}\left\Vert Z\right\Vert
_{L^{2}\left( \mathbb{P};B\right) }\leq c\left( p-1\right) ^{i/2}\left\Vert
Z\right\Vert _{L^{p}\left( \mathbb{P};B\right) }  \label{ZiVsZestimate}
\end{equation}

By H\"{o}lder's inequality for finite sums, 
\begin{equation*}
|Z|_{B}^{q}\leq n^{q-1}\sum_{i=0}^{n}|Z_{i}|_{B}^{q},
\end{equation*}%
and after taking expecations, 
\begin{equation*}
\left\Vert Z\right\Vert _{L^{q}\left( \mathbb{P};B\right) }\leq
n^{1-1/q}\left( \sum_{i=0}^{n}\left\Vert Z_{i}\right\Vert _{L^{q}\left( 
\mathbb{P};B\right) }^{q}\right) ^{1/q}.
\end{equation*}%
Hence, from (\ref{hypercontractivity}) and (\ref{ZiVsZestimate}), 
\begin{eqnarray*}
\left\Vert Z\right\Vert _{L^{q}\left( \mathbb{P};B\right) } &\leq
&cn^{1-1/q}\left( \sum_{i=0}^{n}\left( \frac{q-1}{p-1}\right) ^{\frac{qi}{2}%
}\left\Vert Z_{i}\right\Vert _{L^{p}\left( \mathbb{P};B\right) }^{q}\right)
^{1/q} \\
&\leq &cn^{1-1/q}\left( \sum_{i=0}^{n}\left( q-1\right) ^{\frac{qn}{2}%
}\left\Vert Z\right\Vert _{L^{p}\left( \mathbb{P};B\right) }^{q}\right)
^{1/q} \\
&\leq &cn\left( q-1\right) ^{\frac{n}{2}}\left\Vert Z\right\Vert
_{L^{p}\left( \mathbb{P};B\right) }.
\end{eqnarray*}
\end{proof}

\begin{corollary}
\label{homgeneousnormwienerchaos}Let $X,Y$ be two continuous $G^{N}\left( 
\mathbb{R}^{d}\right) $-valued processes such that, for each $i=1,...,N$,
the projection%
\begin{equation*}
t\mapsto \pi _{i}\left( X_{t}\right) ,\pi _{i}\left( Y_{t}\right) \in \left( 
\mathbb{R}^{d}\right) ^{\otimes i}
\end{equation*}%
belongs to $C_{i}(B_{i})=\mathcal{C}_{0}(B_{i})\oplus \mathcal{C}%
_{1}(B_{i})\oplus ...\mathcal{C}_{n}(B_{i})$, where%
\begin{equation*}
B_{i}=C_{0}\left( \left[ 0,1\right] ,\left( \mathbb{R}^{d}\right) ^{\otimes
i}\right)
\end{equation*}%
and $\mathcal{C}_{j}(B_{i})$ is the $B_{i}$-valued homogeneous Wiener chaos
of degree $j$. Assume $2<q<\infty $. Then there exists a constant $M_{N}$
such that 
\begin{equation*}
\left\Vert d_{\infty }\left( X,Y\right) \right\Vert _{L^{2N}\left( \mathbb{P}%
\right) }\leq \left\Vert d_{\infty }\left( X,Y\right) \right\Vert
_{L^{qN}\left( \mathbb{P}\right) }\leq \sqrt{q}M_{N}\left\Vert d_{\infty
}\left( X,Y\right) \right\Vert _{L^{2N}\left( \mathbb{P}\right) }.
\end{equation*}
\end{corollary}

\begin{proof}
Again, only the second inequality requires a proof. For $i=1,...,N$ define%
\begin{equation*}
Z_{i}:t\mapsto \pi _{i}\left( X_{t}^{-1}\otimes Y_{t}\right) \in \left( 
\mathbb{R}^{d}\right) ^{\otimes i}.
\end{equation*}%
Observe that $Z_{i}\in C_{i}(B_{i})\subset C_{N}\left( B_{i}\right) .$ Let $%
\left\vert .\right\vert _{\infty }$ denote the supremum norm on $\left( 
\mathbb{R}^{d}\right) ^{\otimes i}$ -valued paths. From equation (\ref%
{eqNorms}), 
\begin{equation}
\frac{1}{K_{N}}\max_{i=1,...,N}\left\vert Z_{i}\right\vert _{\infty
}^{1/i}\leq d_{\infty }\left( X,Y\right) \leq
K_{N}\max_{i=1,...,N}\left\vert Z_{i}\right\vert _{\infty }^{1/i}.
\label{eqnorm}
\end{equation}%
Therefore, for all $q\geq 2,$ 
\begin{equation}
\frac{1}{K_{N}}\max_{i=1,...,N}\left\Vert Z_{i}\right\Vert _{L^{qN/i}\left( 
\mathbb{P}\right) }^{1/i}\leq \left\Vert d_{\infty }\left( X,Y\right)
\right\Vert _{L^{qN}\left( \mathbb{P}\right) }\leq
K_{N}\max_{i=1,...,N}\left\Vert Z_{i}\right\Vert _{L^{qN/i}\left( \mathbb{P}%
\right) }^{1/i}.  \label{eqnnorm2}
\end{equation}

From Lemma \ref{LpLqEquivalence} 
\begin{eqnarray*}
\left\Vert Z_{i}\right\Vert _{L^{qN/i}\left( \mathbb{P}\right) }^{1/i} &\leq
&M_{N}^{1/i}\sqrt{qN/i-1}\left\Vert Z_{i}\right\Vert _{L^{2N/i}\left( 
\mathbb{P}\right) }^{1/i} \\
&\leq &\sqrt{q}M_{N}^{\prime }\left\Vert Z_{i}\right\Vert _{L^{2N/i}\left( 
\mathbb{P}\right) }^{1/i} \\
&\leq &\sqrt{q}M_{N}^{\prime \prime }\left\Vert d_{\infty }\left( X,Y\right)
\right\Vert _{L^{2N}\left( \mathbb{P}\right) }
\end{eqnarray*}%
where we used (\ref{eqnnorm2}) with $q=2$ in the last line. Another look at (%
\ref{eqnnorm2}) finishes the proof.
\end{proof}

\section{Large Deviation Results for $\left( \protect\delta _{\protect%
\varepsilon }\mathbf{X}\right) _{\protect\varepsilon >0}$}

From general principles, $\left( \varepsilon X\right) _{\varepsilon >0}$
satisfies a large deviation principle with good rate function $I$ in uniform
topology, where $I$ is given by (see \cite{Ledoux} for example): 
\begin{equation*}
I(y)=\left\{ 
\begin{array}{l}
\frac{1}{2}\left\Vert x\right\Vert _{\mathcal{H}}^{2}\text{ when }y=\emph{i}%
(x)\text{ for some }x\in \mathcal{H} \\ 
+\infty \text{ otherwise.}%
\end{array}%
\right.
\end{equation*}%
It is clear that 
\begin{equation*}
S_{N}\circ \Phi _{m}:\left( C\left( \left[ 0,1\right] ,\mathbb{R}^{d}\right)
,\left\vert \cdot \right\vert _{\infty }\right) \rightarrow C\left( \left[
0,1\right] ,G^{N}\left( \mathbb{R}^{d}\right) ,d_{\infty }\right)
\end{equation*}%
is continuous. By the contraction principle \cite{DZ}, $S_{N}\left(
\varepsilon \Phi _{m}\left( X\right) \right) $ satisfies a large deviation
principle with good rate function 
\begin{equation*}
J_{m}\left( y\right) =\inf \left\{ I\left( x\right) ,\text{ }x\text{ such
that }S_{N}\left( \Phi _{m}\left( x\right) \right) =y\right\} ,
\end{equation*}%
the infimum of the empty set being +$\infty .$ Essentially, a large
deviation principle for $\delta _{\varepsilon }\mathbf{X}$ is obtained by
sending $m$ to infinity. To this end we need

\begin{lemma}
\label{ExpGood}Let $\delta >0$ fixed. Then%
\begin{equation*}
\lim_{m\rightarrow \infty }\overline{\lim_{\varepsilon \rightarrow 0}}%
\varepsilon ^{2}\log \mathbb{P}\left( d_{\infty }\left( S_{N}\left( \Phi
_{m}\left( \varepsilon X\right) \right) ,\delta _{\varepsilon }\mathbf{X}%
\right) >\delta \right) =-\infty .
\end{equation*}
\end{lemma}

\begin{proof}
First observe that 
\begin{equation*}
d_{\infty }\left( S_{N}\left( \Phi _{m}\left( \varepsilon X\right) \right)
,\delta _{\varepsilon }\mathbf{X}\right) =\varepsilon d_{\infty }\left(
S_{N}\left( \Phi _{m}\left( X\right) \right) ,\mathbf{X}\right) .
\end{equation*}%
By standing assumption on $\mathbf{X}$ and lemma \ref{proba_>Lq} we have 
\begin{equation*}
\alpha _{m}:=\left\Vert d_{\infty }\left( S_{N}\left( \Phi _{m}\left(
X\right) \right) ,\mathbf{X}\right) \right\Vert _{L^{2N}}\longrightarrow
_{m\rightarrow \infty }0.
\end{equation*}%
Then, by corollary \ref{homgeneousnormwienerchaos} , 
\begin{equation}
\left\Vert d_{\infty }\left( S_{N}\left( \Phi _{m}\left( X\right) \right) ,%
\mathbf{X}\right) \right\Vert _{L^{q}}\leq M_{N}\sqrt{q}\alpha _{m}\text{ \
\ \ }\forall q\geq 2N.  \label{Lqconv}
\end{equation}%
We then estimate 
\begin{eqnarray*}
\mathbb{P}\left( d_{\infty }\left( S_{N}\left( \Phi _{m}\left( \varepsilon
X\right) \right) ,\delta _{\varepsilon }\mathbf{X}\right) >\delta \right) &=&%
\mathbb{P}\left( d_{\infty }\left( S_{N}\left( \Phi _{m}\left( X\right)
\right) ,\mathbf{X}\right) >\frac{\delta }{\varepsilon }\right) \\
&\leq &\left( \frac{\delta }{\varepsilon }\right) ^{-q}\sqrt{q}^{q}\alpha
_{m}^{q} \\
&\leq &\exp \left[ q\log \left( \frac{\varepsilon }{\delta }\alpha _{m}\sqrt{%
q}\right) \right] ,
\end{eqnarray*}%
and after choosing $q=1/\varepsilon ^{2}$ we obtain, for $\varepsilon $
small enough, 
\begin{equation*}
\varepsilon ^{2}\log \mathbb{P}\left( \sum_{t\in D}d\left( S_{N}\left( \Phi
_{m}\left( \varepsilon X\right) \right) _{t},\delta _{\varepsilon }\mathbf{X}%
_{t}\right) >\delta \right) \leq \log \left( \frac{\alpha _{m}}{\delta }%
\right) .
\end{equation*}%
Now take the limits $\overline{\lim }_{\varepsilon \rightarrow 0}$ and $%
\lim_{m\rightarrow \infty }$ to finish the proof.
\end{proof}

The following theorem is a straight forward application of the extended
contraction principle \cite[Thm 4.2.23]{DZ} and Lemma \ref{ExpGood}. The
point is that although $S_{N}$ is a only a measurable map, defined as a.s.
limit, it has exponentially good approximations given by $\left\{ S_{N}\circ
\Phi _{m}:m\geq 1\right\} $.

\begin{theorem}
\label{Intermediare}Assume that $S_{N}\left( h\right) ,$ defined as the
pointwise limit of $\left( S_{N}\circ \Phi _{m}\right) \left( h\right) _{.}$
when $m\rightarrow \infty $ exists for $h$ such that $I\left( h\right)
<\infty $ (i.e. all $h$ in the Cameron-Martin space), and that for all $%
\Lambda >0,$ 
\begin{equation}
\lim_{m\rightarrow \infty }\sup_{\left\{ h:I\left( h\right) \leq \Lambda
\right\} }d_{\infty }\left[ \left( S_{N}\circ \Phi _{m}\right) \left(
h\right) ,S_{N}\left( h\right) \right] =0.  \label{I_alpha}
\end{equation}%
Then the family $\left( \delta _{\varepsilon }\mathbf{X}\right)
_{\varepsilon >0}$ satisfies a large deviation principle in uniform topology
with good rate function 
\begin{equation*}
J\left( \mathbf{y}\right) =\inf \left\{ I\left( x\right) ,\text{ }x\text{
such that }S_{N}\left( x\right) =\mathbf{y}\right\} .
\end{equation*}
\end{theorem}

The following Corollary allows us to strenghten the topology.

\begin{corollary}
\label{LDPinHolderNorm}Under the above assumptions and the additional
condition 
\begin{equation*}
\exists \alpha ^{\prime }>0,c>0:\mathbb{E}\left( \exp c\left\Vert \mathbf{X}%
\right\Vert _{\alpha ^{\prime }\text{-H\"{o}lder}}^{2}\right) <\infty ,
\end{equation*}%
the family $\left( \delta _{\varepsilon }\mathbf{X}\right) _{\varepsilon >0}$
satisfies a large deviation principle in $\alpha $-H\"{o}lder topology for
all $\alpha \in \lbrack 0,\alpha ^{\prime })$ with good rate function $%
J\left( \mathbf{y}\right) $.
\end{corollary}

\begin{proof}
By the inverse contraction principle \cite[Thm 4.2.4.]{DZ} it suffices to
show that the laws of $\left( \delta _{\varepsilon }\mathbf{X}\right) $ are
exponentially tight in\ $\alpha $-H\"{o}lder topology. But this follows from
the compact embedding of$\ $ 
\begin{equation*}
C^{\alpha ^{\prime }\text{-H\"{o}lder}}\left( \left[ 0,1\right] ,G^{N}\left( 
\mathbb{R}^{d}\right) \right) \subset \subset C^{\alpha \text{-H\"{o}lder}%
}\left( \left[ 0,1\right] ,G^{N}\left( \mathbb{R}^{d}\right) \right) .
\end{equation*}%
and Gauss tails of $\left\Vert \mathbf{X}\right\Vert _{\alpha ^{\prime }%
\text{-H\"{o}lder}}$. Indeed, for every $M>0,$ 
\begin{equation*}
K_{M}=\left\{ x:\left\Vert x\right\Vert _{\alpha ^{\prime }\text{-H\"{o}lder}%
}\leq M\right\}
\end{equation*}%
defines a precompact set w.r.t. $\alpha $-H\"{o}lder topology. Then%
\begin{eqnarray*}
\varepsilon ^{2}\log \left[ \left( \delta _{\varepsilon }\mathbf{X}\right)
_{\ast }\mathbb{P}\right] \left( K_{M}^{c}\right) &=&\varepsilon ^{2}\log 
\mathbb{P}\left[ \left\Vert \delta _{\varepsilon }\mathbf{X}\right\Vert
_{\alpha ^{\prime }\text{-H\"{o}lder }}>M\right] \\
&=&\varepsilon ^{2}\log \mathbb{P}\left[ \left\Vert \mathbf{X}\right\Vert
_{\alpha ^{\prime }\text{-H\"{o}lder}}>M/\varepsilon \right] \\
&\lesssim &-M^{2}.
\end{eqnarray*}
\end{proof}

\begin{remark}
The same proof applies to the (weaker) $p$-variation topology and the
(stronger) modulus topology considered in \cite{FV03}.
\end{remark}

We emphasize the generality of this approach: a large deviation result for
an enhanced Gaussian process holds provided (a) the enhanced Gaussian
process exists and (b) a condition on the reproducing Kernel of the Gaussian
process. \newline
As a warm-up we show how to apply this to enhanced Brownian motion. We will
then discuss a wider class of Gaussian process, which contains fractional
Brownian motion with Hurst parameter $H>\frac{1}{4}.$

\section{LDP for Enhanced Brownian Motion}

Let $X=B$ is be $d$-dimensional Brownian motion. We take $\Phi _{m}$ be the
piecewise linear approximation over at the dyadic dissection of $\left[ 0,1%
\right] $ with mesh $2^{-m}$. Then $S_{2}(\Phi _{m}(B))$ converges uniformly
on $\left[ 0,1\right] $ and in probability to a process $\mathbf{B}$, see 
\cite{LQ,F,FV03}). Large deviations for $\left( \delta _{\varepsilon }%
\mathbf{B}\right) _{\varepsilon }$ were first obtained in in $p$-variation
topology \cite{LQZ}, then in H\"{o}lder and modulus topology \cite{FV03}.
The following proof is considerably quicker. There are only two things to
check.

\begin{enumerate}
\item Gauss tails of $\left\Vert \mathbf{B}\right\Vert _{\alpha \text{-H\"{o}%
lder}}$ for any $\alpha \in \lbrack 0,1/2)$. There are many different ways
to see this. In \cite{FV03} we use the scaling $d\left( \mathbf{B}_{s},%
\mathbf{B}_{t}\right) \overset{\mathcal{D}}{=}\left\vert t-s\right\vert
^{1/2}$ $\left\Vert \mathbf{B}_{0,1}\right\Vert $, noting that $\left\Vert 
\mathbf{B}_{0,1}\right\Vert $ has Gauss tails, followed by an application of
the Garsia-Rodemich-Rumsey lemma. One could prove the result directly by
extending lemma \ref{homgeneousnormwienerchaos} to $\alpha $-H\"{o}lder
metric.

\item Condition (\ref{I_alpha}). But this is a simple consequence of the
fact that the $1$-variation norm of $h$ is bounded by $\sqrt{2I\left(
h\right) }$ and that the $1$-variation of $\Phi _{m}\left( h\right) $ is
bounded by the $1$-variation of $h$. From Theorem \ref{Intermediare} and
Corollary \ref{LDPinHolderNorm} we deduce
\end{enumerate}

\begin{theorem}
Let $\mathbf{B=}$ $\mathbf{B}\left( \omega \right) \in C\left( \left[ 0,1%
\right] ,G^{2}\left( \mathbb{R}^{d}\right) \right) $ denote enhanced
Brownian Motion. Let $\alpha \in \lbrack 0,1/2)$. Then $\left( \delta
_{\varepsilon }\mathbf{B}\right) _{\varepsilon >0}$ satisfies a large
deviation principle in $\alpha $-H\"{o}lder topology with good rate function 
\begin{equation*}
J\left( \mathbf{y}\right) =\inf \left\{ I\left( x\right) ,\text{ }x\text{
such that }S_{N}\left( x\right) =\mathbf{y}\right\} .
\end{equation*}
\end{theorem}

When restricting $\alpha ~$to $\alpha \in \left( 1/3,1/2\right) $ the
universal limit theorem gives the\ Freidlin-Wentzell theorem (in $\alpha $-H%
\"{o}lder topology) as corollary. See \cite{LQZ,FV03,CFV} for details on
this by now classical application of Lyons' limit theorem.

\begin{remark}
The case of enhanced fractional Brownian motion, $H>1/4$, is similar. First,
Gauss tails of the $\alpha $-H\"{o}lder norm, for any $\alpha \in \lbrack
0,H)$, follow from scaling just as above. (This has already been used in 
\cite{FV04, CFV}). Secondly, condition (\ref{I_alpha}) follows from the fact
that the Cameron-Martin space for fBM is embedded in the space of finite $q$%
-variation for any $q$ less than $1/\left( 1/2+H\right) $, see \cite%
{FVcameron}. In any case, enhanced fractional Brownian motion is a special
case of the processes considered in the next section.
\end{remark}

\section{LDP for a Class of Gaussian Processes}

We fix a parameter $H>0.$ We now specialize to zero-mean $\mathbb{R}^{d}$%
-valued Gaussian processes $\{X_{t}:t\in \lbrack 0,1]\}$ with independent
components $x=X^{i},i=1,...d,$ satisfying%
\begin{eqnarray}
E\left( \left\vert x_{s,t}\right\vert ^{2}\right) &\leq &c\left\vert
t-s\right\vert ^{2H},\text{ for all }s<t  \label{ass1} \\
\left\vert E\left( x_{s,s+h}x_{t,t+h}\right) \right\vert &\leq &c\left\vert
t-s\right\vert ^{2H}\left\vert \frac{h}{t-s}\right\vert ^{2},\text{ for all }%
s,t,h\text{ with }h<t-s  \label{ass2}
\end{eqnarray}

for some constant $c>0$. Note that such a process has a.s. $1/p$-H\"{o}lder
continuous sample paths when $p>1/H$. A particular example of such a
Gaussian process is fractional Brownian motion of Hurst parameter $H$. This
condition is the one which allowed Coutin and Qian \cite{CQ} to construct
the lift of $X$ to a rough path. As before $\mathcal{H}$ denotes the
Cameron-Martin space. We take $\Phi _{m}$ to be the piecewise linear
approximation over at the dyadic dissections of $\left[ 0,1\right] $ with
mesh $2^{-m}$.

\begin{definition}
\label{lift}For $H>1/4,$ we define, according to \cite{CQ}, the enhanced
Gaussian process $\mathbf{X}$ to be the (a.s. and in $L^{q}$ for all $q$)
limit in $1/p$-H\"{o}lder norm of $S_{\left[ p\right] }\left( \Phi
_{m}\left( X\right) \right) $, where $1/p<H$.\footnote{%
The convergence in \label{CQ} was proven in $p$-variation topology, but it
is obvious when reading the proof that it also holds in H\"{o}lder norm.}
\end{definition}

To obtain the large deviation principle for $\left( \delta _{\varepsilon }%
\mathbf{X}\right) _{\varepsilon >0}$ there are, as in the case of EBM, only
two things to check.

\begin{enumerate}
\item Gauss tails of $\left\Vert \mathbf{X}\right\Vert _{1/p\text{-H\"{o}lder%
}}$ for $p>1/H,$ that is 
\begin{equation}
\mathbb{E}\exp \left( \alpha \left\Vert \mathbf{X}\right\Vert _{1/p\text{-H%
\"{o}lder}}^{2}\right) <\infty .  \label{gausstail}
\end{equation}%
First note that for all $q\geq 1,$%
\begin{equation*}
\sup_{0\leq s<t\leq 1}\mathbb{E}\left[ \left( \frac{\left\Vert \mathbf{X}%
_{s,t}\right\Vert }{\left\vert t-s\right\vert ^{H}}\right) ^{q}\right]
<\infty .
\end{equation*}
Similar to Corollary \ref{homgeneousnormwienerchaos}, for all $s,t,$ and $%
k\geq 1,$%
\begin{equation*}
\mathbb{E}\left( \left\vert \frac{\left\Vert \mathbf{X}_{s,t}\right\Vert }{%
\left\vert t-s\right\vert ^{H}}\right\vert ^{2k}\right) ^{1/2k}\leq C\sqrt{k}%
E\left( \left\vert \frac{\left\Vert \mathbf{X}_{s,t}\right\Vert }{\left\vert
t-s\right\vert ^{H}}\right\vert ^{q}\right) ^{1/q},
\end{equation*}%
where $q$ is any fixed real greater than $2\left[ 1/H\right] .$ Now expand
the exponential in a power series to see that there exists $\alpha >0$ s.t. $%
\mathbb{E}\exp \left( \alpha \left\Vert \mathbf{X}_{s,t}\right\Vert
^{2}/\left\vert t-s\right\vert ^{2H}\right) <\infty $ for all $s,t$. The
proof of \ref{gausstail} is then finished using the Garsia-Rodemich-Rumsey
lemma, just as in \cite{FV03}.

\item Condition (\ref{I_alpha}). We show (a)\ that for some $q\in \lbrack
1,2)$ the $q$-variation norm of a Cameron-Martin path is bounded by a
constant times its Cameron-Martin norm and (b) a uniform modulus of
continuity for elements in $\left\{ h:I\left( h\right) \leq \Lambda \right\} 
$. A soft argument, spelled out in Theorem \ref{ThmSoft} below, will then
imply (\ref{I_alpha}).
\end{enumerate}

We start by proving that the quadratic variation of Cameron-Martin paths
converges exponentially fast to $0$.

\begin{lemma}
Let $H\in \left( 1/4,1\right) $. Then there exists a constant $C_{H,c}$ such
that for all $n\geq 0$ and all $h\in \mathcal{H}$ 
\begin{equation*}
\sum_{k=0}^{2^{n}-1}\left\vert h_{\frac{k}{2^{n}},\frac{k+1}{2^{n}}%
}\right\vert ^{2}\leq C_{H,c}\left\vert h\right\vert _{\mathcal{H}}^{2}2^{%
\frac{n\left( 1-4H\right) }{2}}
\end{equation*}
\end{lemma}

\begin{proof}
Without loss of generalities, we assume $d=1$. Write $R\left( s,t\right) =%
\mathbb{E}\left( X_{s}X_{t}\right) $ for the covariance. The Cameron Martin
space is the closure of functions of form $h=\sum_{i}a_{i}R(t_{i},\cdot )$
with respect to the norm 
\begin{equation*}
\left\vert h\right\vert _{\mathcal{H}}=\left\vert \sum_{i}a_{i}R(t_{i},\cdot
)\right\vert _{\mathcal{H}}^{2}=\sum_{i,j}a_{i}a_{j}R\left(
t_{i},t_{j}\right) =\mathbb{E}\left( Z^{2}\right)
\end{equation*}%
where $Z$ is the zero-mean Gaussian r.v. given by $Z=\sum_{i}a_{i}X_{t_{i}}$%
. Observe that $h_{s,t}=\mathbb{E}\left( ZX_{s,t}\right) $. We keep $n$
fixed for the rest of the proof.%
\begin{equation*}
Q_{n}:=\sum_{k=0}^{2^{n}-1}\left\vert h_{\frac{k}{2^{n}},\frac{k+1}{2^{n}}%
}\right\vert ^{2}=\sum_{k=0}^{2^{n}-1}\mathbb{E}\left( X_{\frac{k}{2^{n}},%
\frac{k+1}{2^{n}}}Z\right) ^{2}.
\end{equation*}%
Note that if $U,V$ are two zero-mean Gaussian random variables, then%
\begin{equation}
\mathbb{E}\left( U^{2}V^{2}\right) =2\mathbb{E}\left( UV\right) ^{2}+\mathbb{%
E}(U^{2})\mathbb{E}\left( V^{2}\right) \text{ or }Cov\left(
U^{2},V^{2}\right) =2\mathbb{E}\left( UV\right) ^{2}.  \label{gaussian}
\end{equation}

This allows us to write%
\begin{eqnarray*}
Q_{n} &=&\frac{1}{2}\sum_{k=0}^{2^{n}-1}\mathbb{E}\left( Z^{2}\left[ X_{%
\frac{k}{2^{n}},\frac{k+1}{2^{n}}}^{2}-\mathbb{E}\left( X_{\frac{k}{2^{n}},%
\frac{k+1}{2^{n}}}^{2}\right) \right] \right) \\
&=&\frac{1}{2}\mathbb{E}\left( Z^{2}\sum_{k=0}^{2^{n}-1}\left[ X_{\frac{k}{%
2^{n}},\frac{k+1}{2^{n}}}^{2}-\mathbb{E}\left( X_{\frac{k}{2^{n}},\frac{k+1}{%
2^{n}}}^{2}\right) \right] \right) \\
&\leq &\frac{1}{2}\mathbb{E}\left( Z^{4}\right) ^{1/2}\left(
\sum_{k,l=0}^{2^{n}-1}Cov\left( X_{\frac{k}{2^{n}},\frac{k+1}{2^{n}}}^{2},X_{%
\frac{l}{2^{n}},\frac{l+1}{2^{n}}}^{2}\right) \right) ^{1/2} \\
&=&\sqrt{\frac{3}{2}}\mathbb{E}\left( Z^{2}\right) \left(
\sum_{k,l=0}^{2^{n}-1}\mathbb{E}\left( X_{\frac{k}{2^{n}},\frac{k+1}{2^{n}}%
}X_{\frac{l}{2^{n}},\frac{l+1}{2^{n}}}\right) ^{2}\right) ^{1/2}
\end{eqnarray*}%
where we used (\ref{gaussian}) again. The double sum in line just above
splits naturally in a diagonal part, where $k=l$, and twice the sum over $%
k<l $. \ With (\ref{ass1}) we control the diagonal part,%
\begin{equation*}
\sum_{k=0}^{2^{n}-1}\mathbb{E}\left( X_{\frac{k}{2^{n}},\frac{k+1}{2^{n}}%
}^{2}\right) ^{2}\leq c^{2}2^{n\left( 1-4H\right) }.
\end{equation*}

The other part is estimated by using (\ref{ass2}). Indeed,%
\begin{eqnarray*}
\sum_{0\leq k<l\leq 2^{n}-1}\mathbb{E}\left( X_{\frac{k}{2^{n}},\frac{k+1}{%
2^{n}}}X_{\frac{l}{2^{n}},\frac{l+1}{2^{n}}}\right) ^{2}
&=&\sum_{k=0}^{2^{n}-1}\sum_{m=1}^{2^{n}-k-1}\mathbb{E}\left( X_{\frac{k}{%
2^{n}},\frac{k+1}{2^{n}}}X_{\frac{k+m}{2^{n}},\frac{k+m+1}{2^{n}}}\right)
^{2} \\
&\leq &c^{2}\sum_{k=0}^{2^{n}-1}\sum_{m=1}^{2^{n}-k-1}\left\vert \frac{m}{%
2^{n}}\right\vert ^{4H-4}2^{-4n} \\
&\leq &\left( \sum_{m=1}^{\infty }m^{4H-4}\right) c^{2}2^{n\left(
1-4H\right) }.
\end{eqnarray*}%
Putting everything together yields a constanst $C=C\left( H,c\right) $ such
that 
\begin{equation*}
Q_{n}\leq C_{H}\mathbb{E}\left( Z^{2}\right) 2^{\frac{n\left( 1-4H\right) }{2%
}}=C_{H,c}\left\vert h\right\vert _{\mathcal{H}}2^{\frac{n\left( 1-4H\right) 
}{2}}.
\end{equation*}%
This finishes the proof.
\end{proof}

\begin{corollary}
Let $H\in (1/4,1/2]$ and $q\in \left( 1\vee 1/\left( H+1/4\right) ,2\right) $%
. Then $\mathcal{H}$ is continuously embedded in $C_{0}^{q\text{-var}}\left( %
\left[ 0,1\right] ,\mathbb{R}^{d}\right) $. More precisely, there exists a
constant $C_{H,q,c}$ such that for all $h\in \mathcal{H}$, 
\begin{equation*}
\left\vert h\right\vert _{q\text{-var}}\leq C_{H,q,c}\left\vert h\right\vert
_{\mathcal{H}}.
\end{equation*}
\end{corollary}

\begin{proof}
We use the inequality%
\begin{equation*}
\left\vert h\right\vert _{q\text{-var}}^{q}\leq C_{\gamma
,q}\sum_{n=1}^{\infty }n^{\gamma }\sum_{k=0}^{2^{n}-1}\left\vert h_{\frac{k}{%
2^{n}},\frac{k+1}{2^{n}}}\right\vert ^{q},
\end{equation*}%
see \cite{LQ} for example. From the previous lemma and H\"{o}lder's
inequality,%
\begin{eqnarray*}
\sum_{k=0}^{2^{n}-1}\left\vert h_{\frac{k}{2^{n}},\frac{k+1}{2^{n}}%
}\right\vert ^{q} &\leq &\left( \sum_{k=0}^{2^{n}-1}\left\vert h_{\frac{k}{%
2^{n}},\frac{k+1}{2^{n}}}\right\vert ^{2}\right) ^{q/2}2^{n\left(
1-q/2\right) } \\
&\leq &C\left\vert h\right\vert _{\mathcal{H}}^{q}2^{n\left[ \left(
1-4H\right) \frac{q}{4}+1-\frac{q}{2}\right] }.
\end{eqnarray*}%
Observe 
\begin{equation*}
\left( 1-4H\right) \frac{q}{4}+1-\frac{q}{2}<0\Leftrightarrow q\left(
H+1/4\right) >1.
\end{equation*}%
Hence, for every $q\geq 1$ such that $q>1/\left( H+1/4\right) $, we have 
\begin{equation*}
\left\vert h\right\vert _{q\text{-var}}\leq C_{H,q,c}\left\vert h\right\vert
_{\mathcal{H}}.
\end{equation*}
\end{proof}

This result is not optimal. Cameron-Martin paths associated to usual
Brownian motion are of bounded variation ($\equiv $ finite $1$-variation)
whereas the above corollary applied to $H=1/2$ only gives finite $q$%
-variation for $q>4/3$. In \cite{FVcameron}, we prove that Cameron-Martin
paths associated to fractional Brownian motion are of finite $q$-variation
for any $q$ greater than $1\vee 1/\left( H+1/2\right) $.

\begin{lemma}
Let $H\in \left( 0,1\right) $. Then $\mathcal{H}$ is continuously embedded
in $C_{0}^{H\text{-H\"{o}lder}}\left( \left[ 0,1\right] \right) $. More
precisely, there exists a constant $C_{H,c}$ such that for all $h\in 
\mathcal{H}$, 
\begin{equation*}
\left\vert h\right\vert _{H\text{-H\"{o}lder}}\leq C_{H,c}\left\vert
h\right\vert _{\mathcal{H}}.
\end{equation*}
\end{lemma}

\begin{proof}
We use the notation of the above proofs.\ In particular, $h_{s,t}=\mathbb{E}%
\left( ZX_{s,t}\right) .$ By Cauchy-Schwartz and (\ref{ass1}),%
\begin{equation*}
\left\vert h_{s,t}\right\vert \leq \sqrt{\mathbb{E}\left( Z^{2}\right) }%
\sqrt{\mathbb{E}\left( X_{s,t}^{2}\right) }\leq C\left\vert h\right\vert _{%
\mathcal{H}}\left\vert t-s\right\vert ^{H}.
\end{equation*}
\end{proof}

\begin{theorem}
\label{ThmSoft}For any fixed $N\geq 1,$ for all $\Lambda >0,$ 
\begin{equation*}
\lim_{m\rightarrow \infty }\sup_{\left\{ h:I\left( h\right) \leq \Lambda
\right\} }d_{\infty }\left[ \left( S_{N}\circ \Phi _{m}\right) \left(
h\right) ,S_{N}\left( h\right) \right] =0.
\end{equation*}
\end{theorem}

\begin{proof}
Observe first that for $h$ such that $I\left( h\right) \leq \Lambda $, we
have%
\begin{eqnarray*}
\left\vert \Phi _{m}\left( h\right) -h\right\vert _{\infty } &\leq
&\sup_{\left\vert t-s\right\vert \leq 2^{-m}}\left\vert h_{s,t}\right\vert \\
&\leq &\Lambda ^{1/2}C_{H}2^{-mH}.
\end{eqnarray*}%
Observe also that for $q\geq 1/\left( 1/4+H\right) ,$%
\begin{equation*}
\left\vert \Phi _{m}\left( h\right) \right\vert _{q\text{-var}}\leq
3\left\vert h\right\vert _{q\text{-var}}\leq 3C_{q,h}\Lambda ^{1/2}.
\end{equation*}%
By interpolation, for $q^{\prime }>q,$ we therefore obtain that%
\begin{equation*}
\lim_{m\rightarrow \infty }\sup_{\left\{ h:I\left( h\right) \leq \Lambda
\right\} }\left\vert \Phi _{m}\left( h\right) -h\right\vert _{q^{\prime }%
\text{-var}}=0.
\end{equation*}

Then using theorem 2.2.2 in \cite{Ly} (or just the continuity of the Young
integral), we obtain our theorem.
\end{proof}

We have therefore the following:

\begin{corollary}
Assume that for some $H\in \left( 1/4,1\right) $, the Gaussian process
satisfies condition (\ref{ass1}) and (\ref{ass2}) and let $\mathbf{X}$ be as
in definition \ref{lift}.\ Then, $\left( \delta _{\varepsilon }\mathbf{X}%
\right) $ satisfies a large deviation principle in $1/p$-H\"{o}lder topology
where $p>1/H$ with good rate function%
\begin{equation*}
J\left( \mathbf{y}\right) =\inf \left\{ \frac{1}{2}\left\vert x\right\vert _{%
\mathcal{H}},\text{ }x\text{ such that }S_{\left[ p\right] }\left( x\right) =%
\mathbf{y}\right\} .
\end{equation*}
\end{corollary}

Once again, applying the It\^{o} map to $\delta _{\varepsilon }\mathbf{X},$
we can obtain a generalization of the Freidlin-Wentzell theorem, i.e. a
large deviation principle for solution of stochastic differential equation
driven by the enhanced Gaussian process $\mathbf{X}$.

\begin{acknowledgement}
The authors would like to thank Christer Borell, Ben Garling for helpful
discussions and an anonymous referee for his/her careful reading. The first
author acknowledges partial support from a Leverhulme Research Fellowship.
\end{acknowledgement}

\end{document}